\documentclass{notices}

\usepackage{amsmath, amssymb, amsthm, enumitem, esint, graphicx, float}

\newcommand{\lb}{\linebreak[1]}
\newcommand{\tf}{\mathfrak{t}}
\newcommand{\MM}{\mathcal{M}}
\newcommand{\IR}{\mathbb{R}}
\newcommand{\IF}{\mathbb{F}}

\newcommand{\IZ}{\mathbb{Z}}

\newcommand{\RR}{\mathcal{R}}
\newcommand{\NN}{\mathcal{N}}
\renewcommand{\SS}{\mathcal{S}}
\newcommand{\XX}{\mathcal{X}}
\newcommand{\eps}{\varepsilon}

\DeclareMathOperator*{\PSC}{PSC}
\DeclareMathOperator*{\sing}{sing}
\DeclareMathOperator*{\surg}{surg}
\DeclareMathOperator*{\thick}{thick}
\DeclareMathOperator*{\thin}{thin}
\DeclareMathOperator*{\almostsing}{almost-sing}
\DeclareMathOperator*{\can}{can}
\DeclareMathOperator*{\length}{length}
\DeclareMathOperator*{\loc}{loc}
\DeclareMathOperator*{\Bry}{Bry}

\DeclareMathOperator{\Met}{Met}
\DeclareMathOperator{\Diff}{Diff}

\DeclareMathOperator{\Ric}{Ric}
\DeclareMathOperator{\Rm}{Rm}

\DeclareMathOperator{\Isom}{Isom}

\DeclareMathOperator{\vol}{vol}

\newcommand{\dotcup}{\ensuremath{\mathaccent\cdot\cup}}
\newcommand{\EMPTY}[1]{}

\newtheorem{Theorem}{Theorem}[section]

\newtheorem{Conjecture}[Theorem]{Conjecture}

\theoremstyle{definition}
\newtheorem{Definition}[Theorem]{Definition}

\numberwithin{equation}{section}

\title{Recent developments in Ricci flows}
\author{
  Richard H Bamler
  \affil{
    The author is a professor of mathematics at UC Berkeley. His email address is rbamler@berkeley.edu.
    He thanks Robert Bamler, Paula Burkhardt-Guim, Bennett Chow, Bruce Kleiner, Yi Lai and John Lott for helpful feedback on an early version of the manuscript. 
    }
}
\date{\today}

\begin{document}
\maketitle

\section{Introduction}
The beginning of the new millennium marked an important point in the history of geometric analysis.
In a series of three very condensed papers \cite{Perelman1,Perelman2,Perelman3}, Perelman presented a solution of the Poincar\'e and Geometrization Conjectures.
These conjectures --- the former of which was 100 years old and a Millennium Problem --- were fundamental problems concerning the topology of 3-manifolds, yet their solution employed Ricci flow, a technique introduced by Hamilton in 1982, that combines methods of PDE and Riemannian geometry (see below for more details).
These spectacular applications were far from coincidental as they provided a new perspective on 3-manifold topology using the geometric-analytic language of Ricci flows.

It took several years for the scientific community to digest Perelman's arguments and verify his proof.
This was an exciting time for geometric analysts, since the new techniques generated a large number of interesting questions and potential applications.
Amongst these, perhaps the most intriguing goal was to find further applications of Ricci flow to problems in topology.

%

This goal will form the central theme of this article.
I will first summarize the important results pertaining to Perelman's work.
Next, I will discuss more recent research that stands in the same spirit; first in dimension 3, leading, amongst others, to the resolution of the Generalized Smale Conjecture and, second, work in higher dimensions, where topological applications may be forthcoming.
There are also many other interesting recent applications of Ricci flows, which are more of geometric or analytic nature, and which I won't have space to cover here, unfortunately.
Among these are, for example, work on manifolds with positive isotropic curvature \cite{Brendle_PIC_2018}, smoothing constructions for limit spaces of lower curvature bounds \cite{Simon_Topping_local_moll, Bamler_CR_Wilking_2019,Lai_local_nn_2019}, lower scalar curvature bounds for $C^0$-metrics \cite{Bamler_RF_proof_Gromov, Burkhardt_2019}, as well as work on K\"ahler-Ricci flow.

\section{Riemannian geometry and Ricci flow}
A Riemannian metric $g$ on a smooth manifold $M$ is given by a family of inner products on its tangent spaces $T_p M$, $p \in M$, that can be expressed in local coordinates $(x^1, \ldots, x^n)$ as 
\[ g = \sum_{i,j} g_{ij} \, dx^i \, dx^j, \]
where $g_{ij} = g_{ji}$ are smooth coefficient functions.
A metric allows us to define notions such as angles, lengths  and distances on $M$.
One can show that the metric $g$ near any point $p \in M$ looks Euclidean up to order 2, that is we can find suitable local coordinates around $p$ such that
\[ g_{ij} = \delta_{ij} + O(r^2). \]
The second order terms can be described by an invariant called the \emph{Riemann curvature tensor} $\Rm_{iklj}$, which is independent of the choice of coordinates up to a simple transformation rule.
Various components of this tensor have local geometric interpretations.
A particularly interesting component is called the \emph{Ricci tensor} $\Ric_{ij}$, which arises from tracing the $k,l$ indices of $\Rm_{iklj}$.
This tensor roughly describes the second variation of areas of hypersurfaces under normal deformations.

A Ricci flow on a manifold $M$ is given by a smooth family $g(t)$, $t \in [0,T)$, of Riemannian metrics satisfying the evolution equation
\begin{equation} \label{eq_RF}
 \partial_t g(t) = - 2 \Ric_{g(t)}. 
\end{equation}
Expressed in suitable local coordinates, this equation roughly takes the form of a non-linear heat equation in the coefficient functions:
\begin{equation} \label{eq_dtgij_RF}
 \partial_t g_{ij} = \triangle g_{ij} + \ldots 
\end{equation}
In addition, if we compute the time derivative of the Riemannian curvature tensor $\Rm_{g(t)}$, we obtain
\begin{equation} \label{eq_dtRm_RF}
 \partial_t \Rm_{g(t)} = \triangle \Rm_{g(t)} + Q(\Rm_{g(t)}),
\end{equation}
where the last term denotes a quadratic term; its exact form will not be important in the sequel.
Equations (\ref{eq_dtgij_RF}), (\ref{eq_dtRm_RF}) suggest that the metric $g(t)$ becomes smoother or more homogeneous as time moves on, similar to the evolution of a heat equation, under which heat is distributed more evenly across its domain.
On the other hand, the last term in (\ref{eq_dtRm_RF}) seems to indicate that --- possibly at larger scales or in regions of large curvature --- this diffusion property may be outweighed by some other non-linear effects, which could lead to singularities.

But we are getting ahead of ourselves.
Let us first state the following existence and uniqueness theorem, which was established by Hamilton \cite{Hamilton_3_manifolds} in the same paper in which he introduced the Ricci flow equation:

\begin{Theorem}
Suppose that $M$ is compact and let $g_0$ be an arbitrary Riemannian metric on $M$ (called the initial condition). Then:
\begin{enumerate}
\item The evolution equation (\ref{eq_RF}) combined with the initial condition $g(0) = g_0$ has a unique solution $(g(t))_{t \in [0,T)}$ for some maximal $T \in (0, \infty]$. 
\item If $T < \infty$, then we say that the flow \emph{develops a singularity at time $T$} and the curvature blows up:
$\max_M |{\Rm}_{g(t)}| \xrightarrow[t \nearrow T]{} \infty $.
\end{enumerate}
\end{Theorem}

The most basic examples of Ricci flows are those in which $g_0$ is Einstein, i.e. $\Ric_{g_0} = \lambda g_0$.
In this case the flow simply evolves by rescaling:
\begin{equation} \label{eq_homotheticsol}
 g(t) = (1-2\lambda t) g_0
\end{equation}
So for example, a round sphere ($\lambda > 0$) shrinks under the flow and develops a singularity in finite time, where the diameter goes to 0.
On the other hand, if we start with a hyperbolic metric ($\lambda < 0$), then the flow is immortal, meaning that it exists for all times, and the metric expands linearly.
%
In the following we will consider more general initial metrics $g_0$ and hope that --- at least in some cases --- the flow is asymptotic to a solution of the form  (\ref{eq_homotheticsol}).
This will then allow us to understand the topology of the underlying manifold in terms of the limiting geometry.


\section{Dimension 2}
In dimension 2, Ricci flows are very well behaved:

\begin{Theorem} \label{Thm_RF_2d}
Any Ricci flow on a compact 2-dimensional manifold converges, modulo rescaling, to a metric of constant curvature.
\end{Theorem}

In addition, one can show that the flow in dimension 2 preserves the conformal class, i.e. for all times $t$ we have $g(t) = f(t) g_0$ for some smooth positive function $f(t)$ on $M$.
This observation, combined with Theorem~\ref{Thm_RF_2d} can in fact be used to prove of the Uniformization Theorem:\footnote{The proof of Theorem~\ref{Thm_RF_2d}, which was established by Chow and Hamilton \cite{Hamilton_RF_surfaces, Chow_RF_2sphere}, relied on the Uniformization Theorem.
This dependence was later removed by Chen, Lu, Tian \cite{Chen_Lu_Tian_2006}.
}

\begin{Theorem}
Each compact surface $M$ admits a metric of constant curvature in each conformal class.
\end{Theorem}

This is our first topological consequence of Ricci flow.
Of course, the Uniformization Theorem has already been known for about 100 years.
So in order to obtain any \emph{new} topological results, we will need to study the flow in higher dimensions.

\begin{figure*}[!ht]
\centering
\includegraphics[width=118mm]{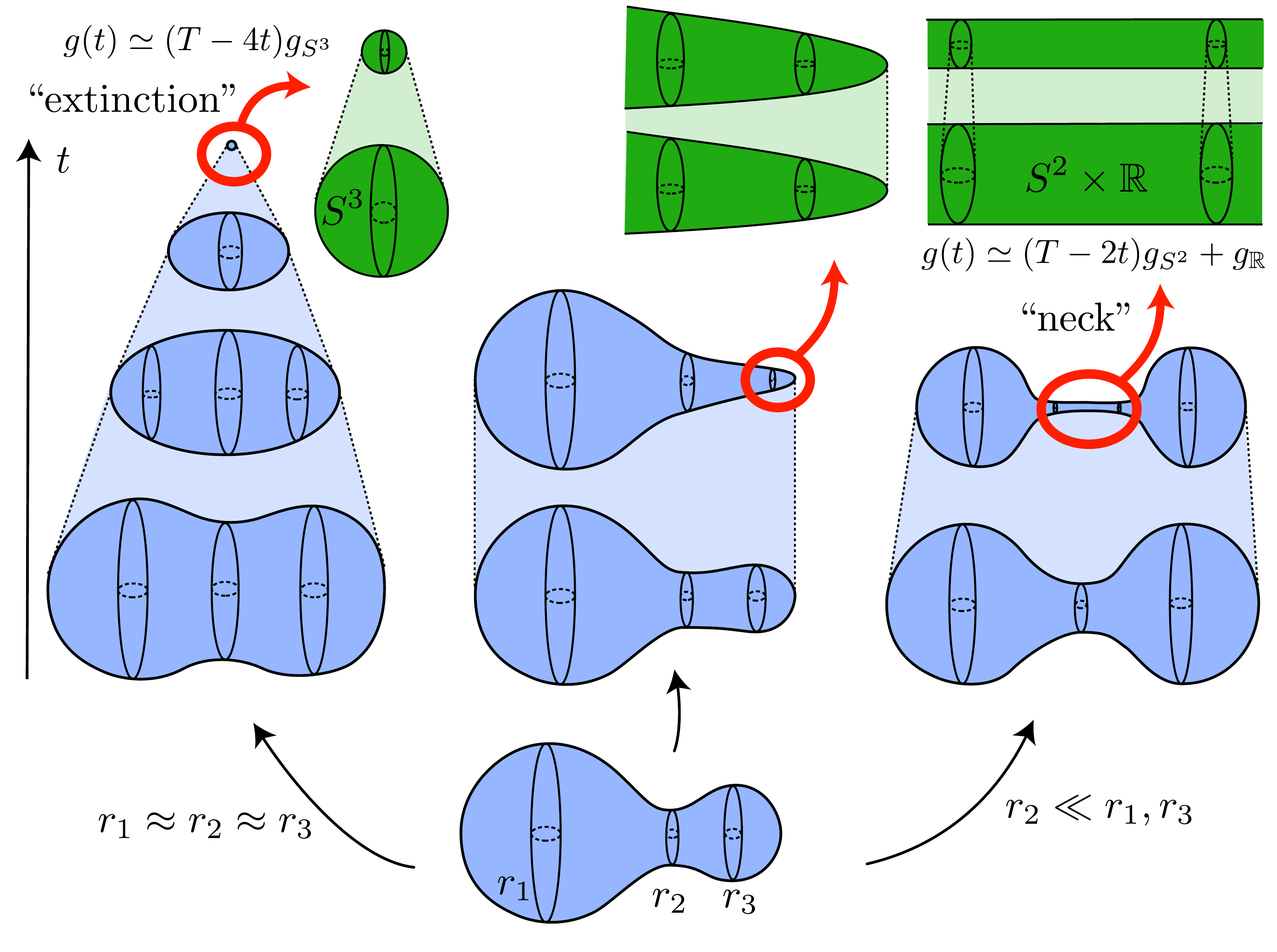}
\caption{Different singularity formations in the rotationally symmetric case, depending on the choice of the radii $r_1, r_2, r_3$. 
The flows depicted on the top are the corresponding singularity models.
These turn out to be the only singularity models, even in the non-rotationally symmetric case (see Subsection~\ref{subsec_sing_models}).
\label{fig_rot_sym}}
\end{figure*}
\section{Dimension 3}
In dimension 3, the behavior of the flow --- and its singularity formation --- becomes far more complicated.
In this section, I will first discuss an example and  briefly review Perelman's analysis of singularity formation and the construction of Ricci flows with surgery.
I will keep this part short and only focus on aspects that will become important later; for a more in-depth discussion and a more complete list of references see the earlier Notices article \cite{Anderson_Notices_04}.
Next, I will focus on more recent work by Kleiner, Lott and the author on singular Ricci flows and their uniqueness and continuous dependence, which led to the resolution of several longstanding topological conjectures.

\subsection{Singularity formation --- an example} \label{subsec_dumbbell} \label{subsec_dumbell}
To get an idea of possible singularity formation of 3-dimensional Ricci flows, it is useful to consider the famous dumbbell example (see Figure~\ref{fig_rot_sym}).
In this example, the initial manifold $(M, g_0)$ is the result of connecting two round spheres of radii $r_1, r_3$ by a certain type of rotationally symmetric neck of radius $r_2$ (see Figure~\ref{fig_rot_sym}).
So $M \approx S^3$ and $g_0 = f^2(s) g_{S^2} + ds^2$ is a warped product away from the two endpoints.

It can be shown that any flow starting from a metric of this form must develop a singularity in finite time.
The singularity \emph{type,} however, depends on the choice of the radii $r_1, r_2, r_3$.
More specifically, if the radii $r_1, r_2, r_3$ are comparable (Figure~\ref{fig_rot_sym}, left), then the diameter of the manifold converges to zero and, after rescaling, the flow 
becomes asymptotically round --- just as in Theorem~\ref{Thm_RF_2d}.
This case is called \emph{extinction}.
On the other hand, if $r_2 \ll r_1, r_3$ (Figure~\ref{fig_rot_sym}, right), then the flow develops a \emph{neck singularity}, which looks like a round cylinder ($S^2 \times \IR$) at small scales.
Note that in this case the singularity only occurs in a certain region of the manifold, while the metric converges to a smooth metric everywhere else.
Lastly, there is also an intermediate case (Figure~\ref{fig_rot_sym}, center), 
in which the flow develops a singularity that is modeled on the Bryant soliton --- a one-ended paraboloid-like model.

\begin{figure}
\centering
\includegraphics[width=60mm]{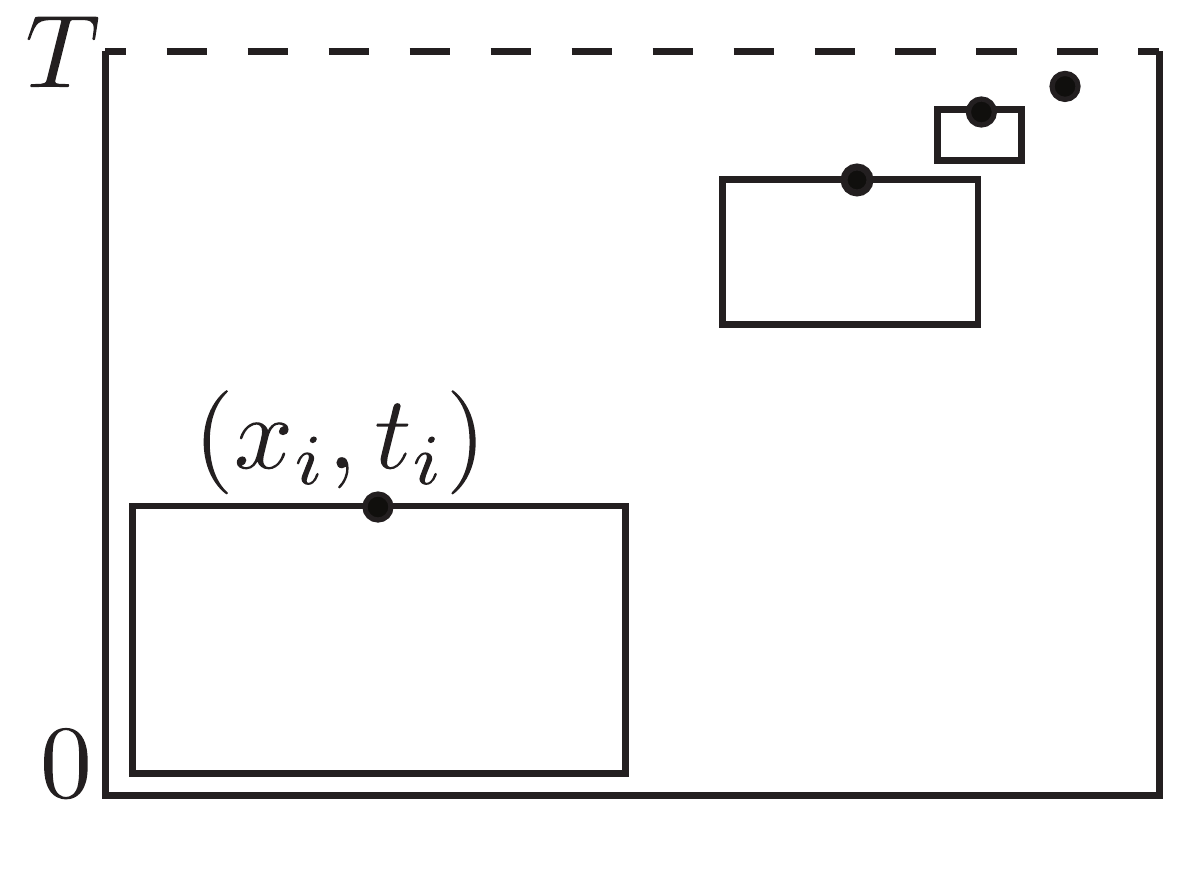}
\caption{A Ricci flow $M \times [0,T)$ that develops a singularity at time $T$ and a sequence of points $(x_i,t_i)$ that ``run into a singularity''.
The geometry in the parabolic neighborhoods around $(x_i,t_i)$ (rectangles) is close to the singularity model modulo rescaling if $i \gg 1$.
\label{fig_blowup}}
\end{figure}
\subsection{Blow-up analysis} \label{subsec_blowup}
Perelman's work implied that the previous example is in fact prototypical for the singularity formation of 
3-dimensional Ricci flows starting from \emph{general,} not necessarily rotationally symmetric, initial data.
In order to make this statement more precise, let us first recall a method called \emph{blow-up analysis,} which is used frequently to study singularities in geometric analysis.

Suppose that $(M, (g(t))_{t \in [0,T)})$ is a Ricci flow  that develops a singularity at time $T < \infty$ (see Figure~\ref{fig_blowup}).
Then we can find a sequence of spacetime points $(x_i, t_i) \in M \times [0,T)$ such that $\lambda_i := |{\Rm}|(x_i, t_i) \to \infty$ and $t_i \nearrow T$; we will say that the points $(x_i, t_i)$ ``run into a singularity''.
Our goal will be to understand the local geometry at small scales near $(x_i, t_i)$, for large $i$.
For this purpose, we 
consider the sequence of pointed, parabolically rescaled flows
\[ \big(M, (g'_i(t) := \lambda_i g(\lambda_i^{-1} t + t_i))_{t \in [-\lambda_i t_i, 0]}, (x_i,0) \big). \]
Geometrically, the flow $(g'_i(t))$ is the result of rescaling distances by $\lambda_i^{1/2}$, times by $\lambda_i$ and an application of a time-shift such that the point $(x_i,0)$ in the new flow corresponds the point $(x_i, t_i)$ in the old flow.
The new flow $(g'_i(t))$ still satisfies the Ricci flow equation and it is defined on larger and larger backwards time-intervals of size $\lambda_i t_i \to \infty$.
Moreover, we have $|{\Rm}|(x_i,0)=1$ on this new flow.
Observe also that the geometry of the original flow near $(x_i,t_i)$ at scale $\lambda_i^{-1/2} \ll 1$ is a rescaling of the geometry of $(g'_i(t))$ near $(x_i,0)$ at scale $1$.

The hope is now to apply a compactness theorem (\`a la Arzela-Ascoli) such that, after passing to a subsequence, we have convergence
\begin{multline} \label{eq_RF_convergence}
 \big( M, (g'_i(t))_{t \in [-\lambda_i^{2} t_i, 0]}, (x_i,0) \big) \\ \xrightarrow[i \to \infty]{} \big( M_\infty, (g_\infty(t))_{t \leq 0}, (x_\infty, 0) \big). 
\end{multline}
The limit is called a \emph{blow-up} or \emph{singularity model}, as it gives valuable information on the singularity near the points $(x_i, t_i)$.
This model is a Ricci flow that is defined for all times $t \leq 0$; it is therefore called \emph{ancient.}
So in summary, a blow-up analysis reduces the study of singularity \emph{formation} to the classification of ancient singularity \emph{models.}

The notion of convergence in (\ref{eq_RF_convergence}) is a generalization of Cheeger-Gromov convergence to Ricci flows, due to Hamilton.
Instead of demanding global convergence of the metric tensors, as in Theorem~\ref{Thm_RF_2d}, we only require convergence up to diffeomorphisms here.
We roughly require that we have convergence
\begin{equation} \label{eq_Hamilton_phi_Cinfty}
 \phi_i^* g'_i (t) \xrightarrow[i \to \infty]{C^\infty_{\loc}} g_\infty
\end{equation}
on $M_\infty \times (-\infty, 0]$ of the pull backs of $g'_i (t)$ via (time-independent) diffeomorphisms $\phi_i : U_i \to V_i \subset M$ that are defined over larger and larger subsets $U_i \subset M_\infty$ and satisfy $\phi_i(x_\infty) = x_i$.
%
We will see later (in Section~\ref{sec_dim_n_geq_4}) that this notion of convergence is too strong to capture the more subtle singularity formation of higher dimensional Ricci flows and we will discuss necessary refinements.
Luckily, in dimension 3 the current notion is still sufficient for our purposes, though.

\subsection{Singularity models and canonical neighborhoods} \label{subsec_sing_models}
%
Arguably one of the most groundbreaking discoveries of Perelman's work was the classification of singularity models of 3-dimensional Ricci flows and the resulting structural description of the flow near a singularity.
The following theorem\footnote{Perelman proved a version of Theorem~\ref{Thm_classificaton_3d} that contained a more qualitative characterization in Case~3., which was sufficient for most applications.
Later, Brendle \cite{Brendle_Bryant_2020} showed that the only possibility in Case~3. is the Bryant soliton.} summarizes this classification.


%
%
%
%
%
%

\begin{Theorem} \label{Thm_classificaton_3d}
Any singularity model $(M_\infty, (g_\infty(t))_{t \leq 0})$ obtained as in (\ref{eq_RF_convergence}) is isometric, modulo rescaling, to one the following:
\begin{enumerate}
\item A quotient of the round shrinking sphere $(S^3, (1-4t)g_{S^3})$.
\item The round shrinking cylinder $(S^2 \times \IR, (1 - 2t)g_{S^2} + g_{\IR})$ or its quotient $(S^2 \times \IR)/\IZ_2$.
\item The Bryant soliton $(M_{\Bry}, (g_{\Bry}(t)))$
\end{enumerate}
\end{Theorem}

Note that these three models correspond to the three cases in the rotationally symmetric dumbbell example from Subsection~\ref{subsec_dumbbell} (see Figure~\ref{fig_rot_sym}).
The Bryant soliton in 3. is a rotationally symmetric solution to the Ricci flow on $\IR^3$ with the property that all its time-slices are isometric to a metric of the form
\[ g_{\Bry} = f^2(r) g_{S^2} + dr^2, \qquad f(r) \sim \sqrt{r}. \]
The name \emph{soliton} refers to the fact that all time-slices of the flow are isometric, so the flow merely evolves by pullbacks of a family of diffeomorphisms (we will encounter this soliton again at the end of this article).

It turns out that we can do even better than Theorem~\ref{Thm_classificaton_3d}: we can describe the structure of the flow near \emph{any} point of the flow that is close to a singularity --- not just near the points used to construct the blow-ups.
In order to state the next result efficiently, we will need to consider the class of \emph{$\kappa$-solutions.}
This class simply consists of all solutions listed in Theorem~\ref{Thm_classificaton_3d}, plus an additional compact, ellipsoidal solution (the details of this solution won't be important here\footnote{This solution does not occur as a singularity of a single flow, but can be observed as some kind of transitional model in families of flows that interpolate between two different singularity models.
}). 
Then we have:

\begin{Theorem}[Canonical Neighborhood Theorem] \label{Thm_CNT}
If $(M, (g(t))_{t \in [0,T)})$, $T < \infty$, is a 3-dimensional Ricci flow and $\eps >0$, then there is a constant $r_{\can} (g(0), T, \eps) > 0$ such that for any $(x,t) \in M \times [0,T)$ with the property that
\[ r := |{\Rm}|^{-1/2}(x,t) \leq r_{\can} \]
the geometry of the metric $g(t)$ restricted to the ball $B_{g(t)}(x, \eps^{-1} r)$ is $\eps$-close\footnote{Similar to (\ref{eq_Hamilton_phi_Cinfty}), this roughly means that there is a diffeomorphism between an $\eps^{-1}$-ball in a $\kappa$-solution and this ball such that the pullback of $r^{-2} g(t)$ is $\eps$-close in the $C
^{[\eps^{-1}]}$-sense to the metric on the $\kappa$-solution.} to a time-slice of a $\kappa$-solution.
\end{Theorem}

Let us digest the content of this theorem.
Recall that the norm of the curvature tensor $|{\Rm}|$, which can be viewed as a measure of ``how singular the flow is near a point'', has the dimension of $\length^{-2}$, so $r = |{\Rm}|^{-1/2}$ can be viewed as a ``curvature scale''.
So, in other words, Theorem~\ref{Thm_CNT} states that regions of high curvature locally look very much like cylinders, Bryant solitons, round spheres etc.

\begin{figure}
\centering
\includegraphics[width=60mm]{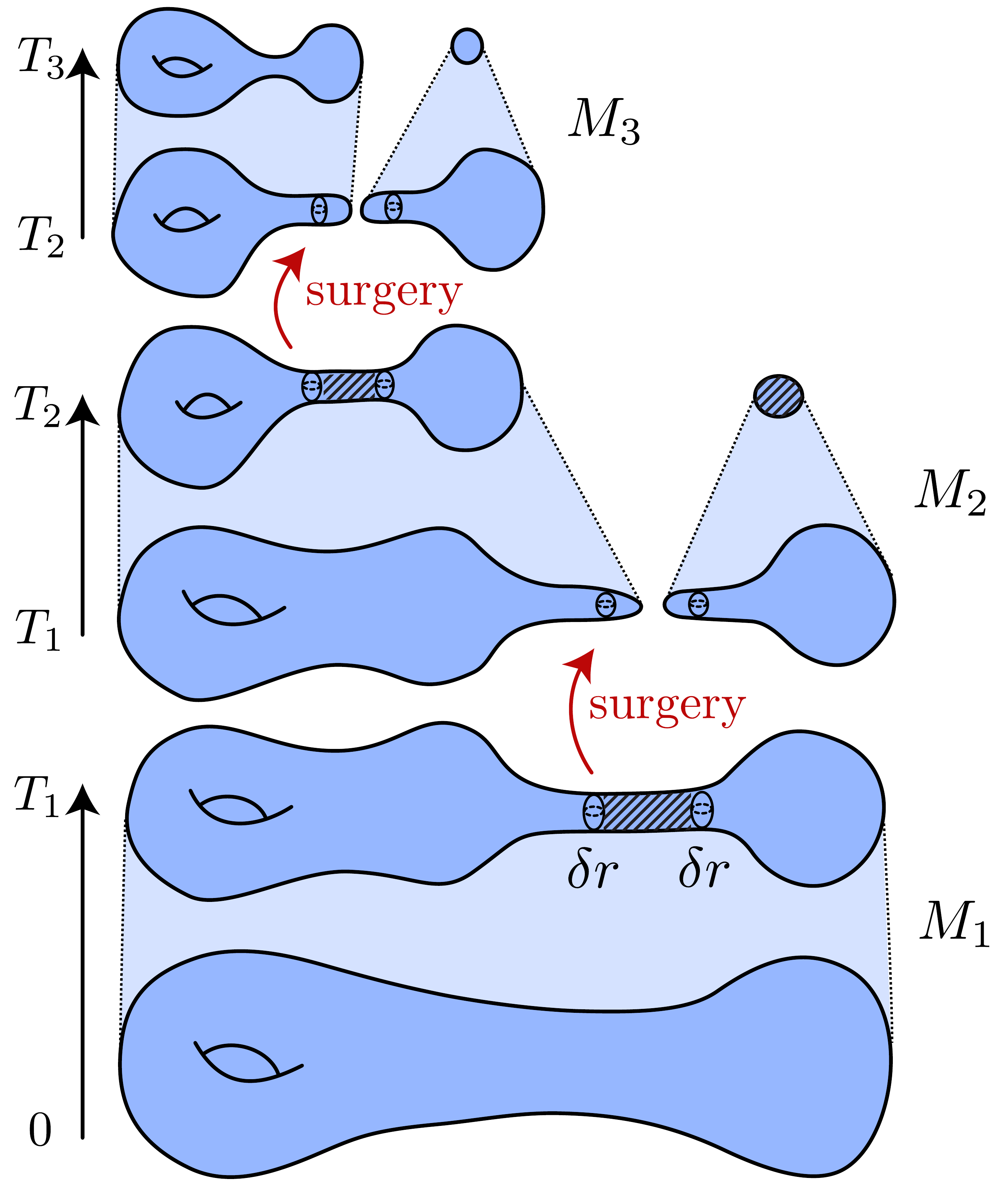}
\caption{A schematic depiction of a Ricci flow with surgery. The almost-singular parts $M_{\text{almost-sing}}$, i.e. the parts that are discarded under each surgery construction, are hatched.\label{fig_RF_surgery}}
\end{figure}
\subsection{Ricci flow with surgery} \label{subsec_RFsurg}
Our understanding of the structure of the flow near a singularity now allows us to carry out a so-called \emph{surgery construction.} 
Under this construction (almost) singularities of the flow are removed, resulting in a ``less singular'' geometry, from which the flow can be restarted.
This will lead to a new type of flow that is defined beyond its singularities and which will provide important information on the underlying manifold.


Let me be more precise.
A \emph{(3-dimensional) Ricci flow with surgery} (see Figure~\ref{fig_RF_surgery}) consists of a sequence of Ricci flows
\begin{align*}
&(M_1, (g_1(t))_{t \in [0,T_1]}), \quad (M_2, (g_2(t))_{t \in [T_1,T_2]}), \quad \\
&\qquad (M_3, (g_3(t))_{t \in [T_2,T_3]}), \quad \ldots, 
\end{align*}
which live on manifolds $M_1, M_2, \ldots$ of possibly different topology and are parameterized by consecutive time-intervals of the form $[0,T_1], [T_1, T_2], \ldots$ whose union equals $[0, \infty)$.
The time-slices $(M_i, g_i(T_i))$ and $(M_{i+1}, g_{i+1}(T_i))$ are related by a \emph{surgery process,} which can be roughly summarized as follows.
Consider the set $M_{\almostsing} \subset M_i$ of all points of high enough curvature, such that they have a canonical neighborhood as in Theorem~\ref{Thm_CNT}.
Cut $M_i$ open along approximate cross-sectional 2-spheres of diameter $r_{\surg}(T_i) \ll 1$ near the cylindrical ends of $M_{\almostsing}$, discard most of the high-curvature components (including the closed, spherical components of $M_{\almostsing}$), and glue in cap-shaped 3-disks to the cutting surfaces. 
%
In doing so we have constructed a new, ``less singular'', Riemannian manifold $(M_{i+1}, g_{i+1}(T_i))$, from which we can restart the flow.
Stop at some time $T_{i+1} > T_i$, shortly before another singularity occurs and repeat the process.

The precise surgery construction is quite technical and more delicate than presented here.
The main difficulty in this construction is to ensure that the surgery times $T_i$ do not accumulate, i.e. that the flow can be extended for all times.
It was shown by Perelman that this and other difficulties can indeed be overcome:

\begin{Theorem}
Let $(M,g)$ be a closed, 3-dimensional Riemannian manifold.
If the surgery scales $r_{\surg}(T_i) > 0$ are chosen sufficiently small (depending on $(M,g)$ and $T_i$), then a Ricci flow with surgery with initial condition $(M_1, g_1(0)) = (M,g)$ can be constructed.
\end{Theorem}

Note that the topology of the underlying manifold $M_i$ may change in the course of a surgery, but only in a controlled way.
In particular, it is possible to show that for any $i$ the initial manifold $M_1$ is diffeomorphic to a connected sum of components of $M_i$ and copies of spherical space forms $S^3/\Gamma$ and $S^2 \times S^1$.
So if the flow goes \emph{extinct} in finite time, meaning that $M_i = \emptyset$ for some large $i$, then
\begin{equation} \label{eq_sph_prime_dec}
 M_1 \approx \#_{j=1}^k (S^3/\Gamma_j) \# m (S^2 \times S^1). 
\end{equation}
Perelman moreover showed that if $M_1$ is simply connected, then the flow \emph{has} to go extinct and therefore $M_1$ must be of the form (\ref{eq_sph_prime_dec}).
This immediately implies the Poincar\'e Conjecture --- our first true topological application of Ricci flow:

%

\begin{Theorem}[Poincar\'e Conjecture]
Any simply connected, closed 3-manifold is diffeomorphic to $S^3$.
\end{Theorem}

On the other hand, Perelman showed that if the Ricci flow with surgery does not go extinct, meaning if it exists for all times, then for large times $t \gg 1$ the flow decomposes the manifold (at time $t$) into a thick and a thin part:
\begin{equation} \label{eq_thick_thin_3d}
 M_{\thick}(t) \,\, \dotcup \,\, M_{\thin}(t), 
\end{equation}
such that the metric on $M_{\thick}(t)$ is asymptotic to a hyperbolic metric and the metric on $M_{\thin}(t)$ is collapsed along fibers\footnote{The precise description of these fibrations is a bit more complex and omitted here. In addition, $M_{\thick}(t)$ and $M_{\thin}(t)$ are separated by embedded 2-tori that are incompressible, meaning that the inclusion into the ambient manifold is an injection on $\pi_1$.
} of type $S^1, S^2$ or $T^2$.
A further topological analysis of these induced fibrations  implied the Geometrization Conjecture --- our second topological application of Ricci flow:

\begin{Theorem}[Geometrization Conjecture]
Every closed 3-manifold is a connected sum of manifolds that can be cut along embedded, incompressible copies of $T^2$ into pieces which each admit a locally homogeneous geometry.
\end{Theorem}

\begin{figure}
\hspace{-3mm}
\includegraphics[width=70mm]{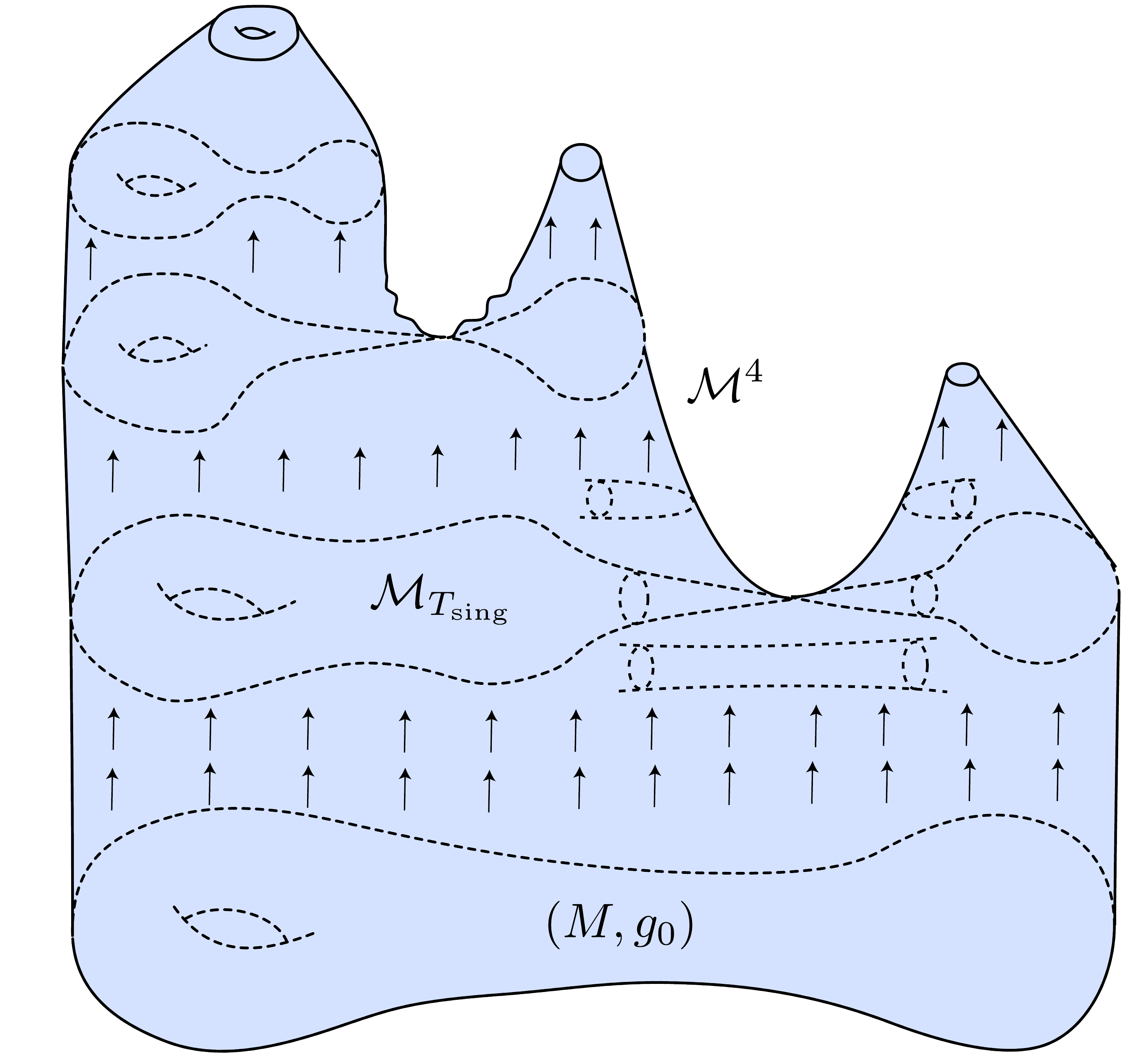}
\caption{Illustration of a singular Ricci flow given by a Ricci flow spacetime.
The arrows indicate the time-vector field $\partial_\tf$.\label{fig_sing_RF}}
\end{figure}
\subsection{Ricci flows through singularities} \label{subsec_RF_through_sing}
Despite their spectacular applications, Ricci flows with surgery have one major drawback: their construction is not canonical.
In other words, each surgery step depends on a number of auxiliary parameters, for which there does not seem to be a canonical choice, such as:
\begin{itemize}
\item The surgery scales $r_{\surg}(T_i)$, i.e. the diameters of the cross-sectional spheres along which the manifold is cut open. 
These scales need to be positive and small.
\item The precise locations of these surgery spheres.
\end{itemize}
Different choices of these parameters may influence the future development of the flow significantly (as well as the space of future surgery parameters). 
Hence a Ricci flow with surgery is not \emph{uniquely} determined by its initial metric.

This disadvantage was already recognized in Perelman's  work, where he conjectured that there should be another flow, in which surgeries are effectively carried out automatically at an infinitesimal scale (think ``$r_{\surg} = 0$''), or which, in other words, ``flows through singularities''.

Perelman's conjecture was recently resolved by Kleiner, Lott and the author \cite{Kleiner_Lott_singular, bamler_kleiner_uniqueness_stability}:\footnote{The ``Existence'' part is due to Kleiner and Lott and the ``Uniqueness'' part is due to Kleiner and the author. Part 2. of Theorem~\ref{Thm_sing_RF} follows from a combination of both papers.}

\begin{Theorem} \label{Thm_sing_RF}
There is a notion of singular Ricci flow (through singularities) such that:
\begin{enumerate}
\item For any compact, 3-dimensional Riemannian manifold $(M,g)$ there is a unique singular Ricci flow $\mathcal{M}$ whose initial time-slice $(\MM_0, g_0)$ is $(M,g)$.
\item Any Ricci flow with surgery starting from $(M,g)$ can be viewed as an approximation of $\mathcal{M}$.
More specifically, if we consider a sequence of Ricci flows with surgery starting from $(M,g)$ with surgery scales $\max_t r_{\surg}(t) \to 0$, then these flows converge to $\mathcal{M}$ in a certain sense.
\end{enumerate}
\end{Theorem}

In addition, the concept of a singular Ricci flow is far less technical than that of a Ricci flow with surgery ---
in fact, I will be able to state its full definition here.
To do this, I will first define the concept of a Ricci flow spacetime.
In short, this is a smooth 4-manifold that locally looks like a Ricci flow, but which may have  non-trivial \emph{global} topology (see Figure~\ref{fig_sing_RF}).

\begin{Definition} \label{Def_RF_spacetime}
A \emph{Ricci flow spacetime} consists of:
\begin{enumerate}
\item A smooth 4-dimensional manifold $\mathcal{M}$ with boundary, called \emph{spacetime.}
\item A \emph{time-function} $\tf : \mathcal{M} \to [0, \infty)$.
Its level sets $\mathcal{M}_t := \mathfrak{t}^{-1}(t)$ are called \emph{time-slices} and we require that $\mathcal{M}_0 = \partial \mathcal{M}$.
\item A \emph{time-vector field} $\partial_{\tf}$ on $\MM$ with $\partial_{\tf} \cdot \tf \equiv 1$.
Trajectories of $\partial_{\tf}$ are called \emph{worldlines.}
\item A family $g$ of inner products on $\ker d\tf \subset T \MM$, which induce a Riemannian metric $g_t$ on each time-slice $\MM_t$.
We require that the Ricci flow equation holds:
\[ \mathcal{L}_{\partial_{\tf}} g_t = - 2\Ric_{g_t}. \]
\end{enumerate}
By abuse of notation, we will often write $\MM$ instead of $(\MM, \tf, \partial_{\tf}, g)$.
\end{Definition}

A classical, 3-dimensional Ricci flow $(M, (g(t))_{t \in [0,T)})$ can be converted into a Ricci flow spacetime by setting $\MM := M \times [0,T)$, letting $\tf, \partial_{\tf}$ be the projection onto the second factor and the pullback of the unit vector field on the second factor, respectively, and letting $g_t$ be the metric corresponding to $g(t)$ on $M \times \{ t \} \approx M$.
Hence worldlines correspond to curves of the form $t \mapsto (x,t)$.

Likewise, a Ricci flow with surgery, given by flows $(M_1,  \lb  (g_1(t))_{t \in [0,T_1]}), \linebreak[0] (M_2,\lb (g_2(t))_{t \in [T_1,T_2]}), \ldots$ can be converted into a Ricci flow spacetime as follows.
Consider first the Ricci flow spacetimes $M_1 \times [0,T_1], M_2 \times [T_1,T_2], \ldots$ arising from each single flow.
We can now glue these flows together by identifying the set of points $U_i^- \subset M_i \times \{ T_i \}$ and $U_i^+ \subset M_{i+1} \times \{ T_i \}$ that survive each surgery step via maps $\phi_i : U_i^- \to U_i^+$.
%
%
%
The resulting space has a boundary that consists of the time-$0$-slice $M_1 \times \{ 0 \}$ and the points
\[ \mathcal{S}_i = (M_i \times \{ T_i \} \setminus U_i^- ) \cup ( M_{i+1} \times \{ T_{i} \} \setminus U_i^+ ),  \]
which were removed and added during each surgery step.
After removing these points, we obtain a Ricci flow spacetime of the form:
\begin{multline} \label{eq_MM_from_RF_surg}
 \MM = (M_1 \times [0,T_1] \cup_{\phi_1} M_2 \times [T_1, T_2] \cup_{\phi_2} \ldots) \\ \setminus (\mathcal{S}_1 \cup \mathcal{S}_2 \cup \ldots). 
\end{multline}
Note that for any regular time $t \in (T_{i-1}, T_i)$ the time-slice $\MM_t$ is isometric to $(M_{i}, g_i(t))$.
On the other hand, the time-slices $\MM_{T_i}$ corresponding to surgery times are incomplete; they have cylindrical open ends of scale $\approx r_{\surg}(T_i)$.

The following definition captures this incompleteness:

\begin{Definition}
A Ricci flow spacetime is \emph{$r$-complete,} for some $r \geq 0$, if the following holds.
Consider a smooth path $\gamma : [0,s_0) \to \MM$ with the property that
\[ \inf_{s \in [0,s_0)} |{\Rm}|^{-1/2} (\gamma(s)) > r \]
and:
\begin{enumerate}
\item $\gamma([0,l)) \subset \MM_t$ is contained in a single time-slice and its length measured with respect to the metric $g_t$ is finite, or
\item $\gamma$ is a worldline, i.e. a trajectory of $\pm\partial_{\tf}$.
\end{enumerate}
Then the limit $\lim_{s \nearrow s_0} \gamma(s)$ exists.
\end{Definition}

So $\MM$ being $r$-complete roughly means that it has only ``holes'' of scale $\lesssim r$.
For example, the flow from (\ref{eq_MM_from_RF_surg}) is $C\max_t r_{surg}(t)$-complete for some universal $C < \infty$.

In addition, Theorem~\ref{Thm_CNT} motivates the following definition:

\begin{Definition}
A Ricci flow spacetime is said to satisfy the \emph{$\epsilon$-canonical neighborhood assumption at scales $(r_1, r_2)$} if for any point $x \in \MM_t$ with $r := |{\Rm}|^{-1/2}(x) \in (r_1, r_2)$ the metric $g_t$ restricted to the ball $B_{g_t} (x, \eps^{-1} r)$ is $\eps$-close, after rescaling by $r^{-2}$, to a time-slice of a $\kappa$-solution.
\end{Definition}

We can finally define singular Ricci flows (through singularities), as used in Theorem~\ref{Thm_sing_RF}:

\begin{Definition}
A \emph{singular Ricci flow} is a Ricci flow spacetime $\MM$ with the following two properties:
\begin{enumerate}
\item It is $0$-complete.
\item For any $\eps  > 0$ and $T < \infty$ there is an $r(\eps, T) > 0$ such that the flow $\MM$ restricted to $[0,T)$ satisfies the $\eps$-canonical neighborhood assumption at scales $(0, r)$.
\end{enumerate}
\end{Definition}

See again Figure~\ref{fig_sing_RF} for a depiction of a singular Ricci flow.
The time-slices $\MM_t$ for $t < T_{\sing}$ develop a cylindrical region, which collapses to some sort of topological double cone singularity in the time-$T_{\sing}$-slice $\MM_{T_{\sing}}$.
This singularity is immediately resolved and the flow is smooth for all $t > T_{\sing}$.

Let us digest the definition of a singular Ricci flow a bit more.
It is tempting to think of the time function $\tf$ as a Morse function and compare critical points with infinitesimal surgeries.
However, this comparison is flawed:
First, by definition $\tf$ cannot have critical points since $\partial_\tf \, \tf = 1$.
In fact, a singular Ricci flow is a completely smooth object.
The ``singular points'' of the flow are not part of $\MM$, but can be obtained after metrically completing each time-slice by adding a discrete set of points.
Second, it is currently unknown whether the set of singular times, i.e. the set of times whose time-slices are incomplete, is discrete.

Similar notions of singular flows have been developed for the mean curvature flow (a close cousin of the Ricci flow).
These are called level set flows and Brakke flows.
However, their definitions differ from singular Ricci flows in that they characterize the flow equation at singular points via barrier and weak integral conditions, respectively.
This is possible, in part, because a mean curvature flow is an embedded object and its singular set has an analytic meaning.
By contrast, the definition of a singular Ricci flow only characterizes the flow on its regular part.
In lieu of a weak formulation of the Ricci flow equation on the singular set, we have to impose the canonical neighborhood assumption, which serves as an asymptotic characterization near the incomplete ends.

Finally, I will briefly explain how singular Ricci flows are constructed and convey the meaning of Part 2. of Theorem~\ref{Thm_sing_RF}.
%
Fix an initial time-slice $(M,g)$ and consider a sequence of Ricci flow spacetimes $\MM^j$ that correspond to Ricci flows with surgery starting from $(M,g)$, with surgery scale $\max_t r_{\surg,j}(t) \to 0$.
It can be shown that these flows are $C \max_t r_{\surg,j}(t)$-complete and satisfy the $\eps$-canonical neighborhood assumption at scales $(C_\eps \max_t r_{\surg,j}(t), r_\eps)$, where $C, C_\eps, r_{\eps}$ do not depend on $j$.
A compactness theorem implies that a subsequence of the spacetimes $\MM^j$ converges to a spacetime $\MM$, which is a singular Ricci flow.
This implies the existence of $\MM$; the proof of uniqueness uses other techniques, which are outside the scope of this article.

\begin{figure*}
\centering
\includegraphics[width=120mm]{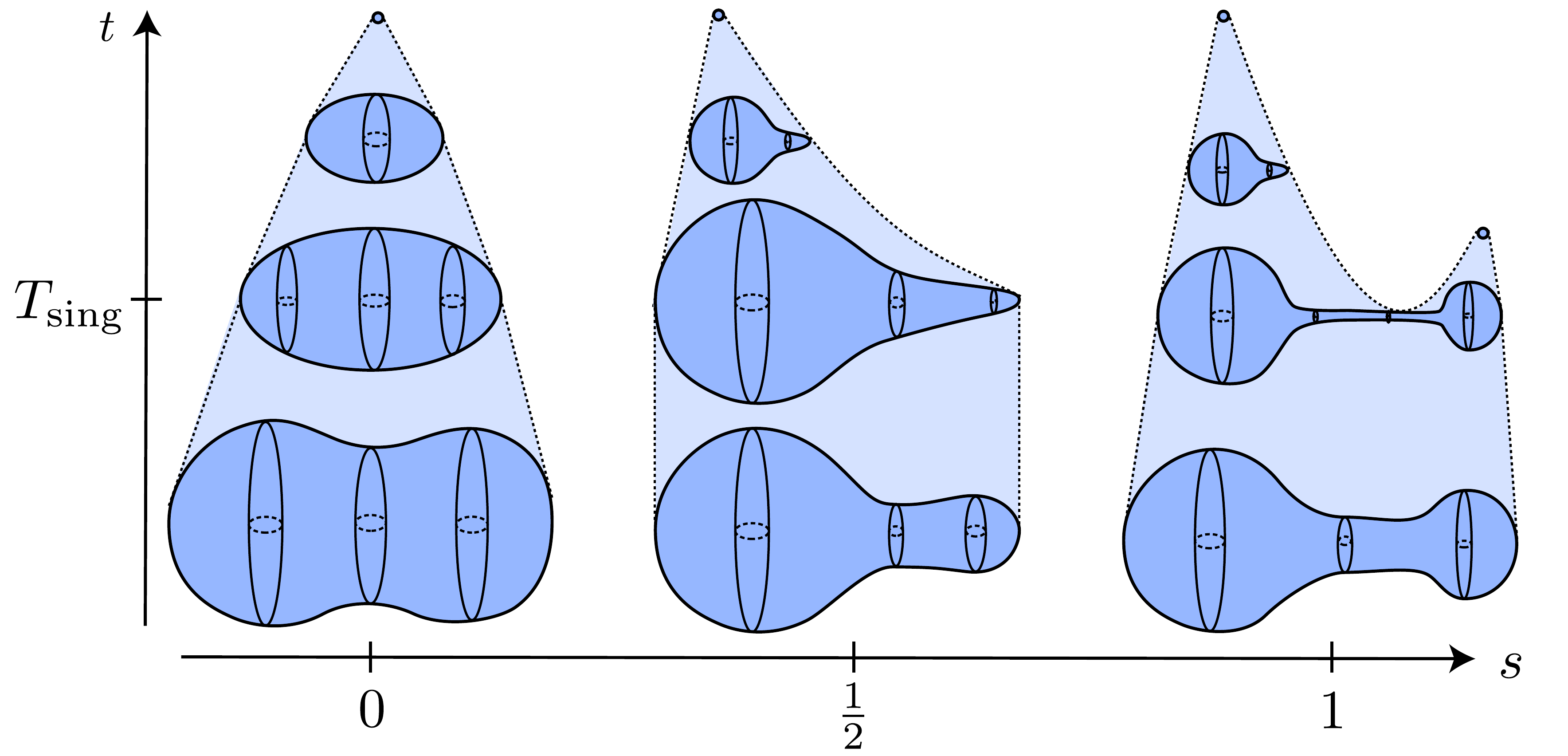}
\caption{A family of singular Ricci flows starting from a continuous family of initial conditions.\label{fig_rot_sym_family}}
\end{figure*}
\subsection{Continuous dependence}
The proof of the uniqueness property in Theorem~\ref{Thm_sing_RF}, due to Kleiner and the author, implies an important continuity property, which will lead to further topological applications.
To state this property, let $M$ be a compact 3-manifold and for every Riemannian metric $g$ on $M$ let $\MM^g$ be the singular Ricci flow with initial condition $(\MM^g_0, g) = (M,g)$.

\begin{Theorem} \label{Thm_cont_dep}
The flow $\MM^g$ depends continuously on $g$.
\end{Theorem}

Recall that the topology of the flow $\MM^g$ may change as we vary $g$.
We therefore have to choose an appropriate sense of continuity in Theorem~\ref{Thm_cont_dep} that allows such a topological change.
This is roughly done via a topology and lamination structure on the disjoint union $\bigsqcup_g \MM^g$, transverse to which the variation of the flow can be studied locally.

Instead of diving into these technicalities, let us discuss the example illustrated in Figure~\ref{fig_rot_sym_family}.
In this example $(g_s)_{s \in [0,1]}$ denotes a continuous family of metrics on $S^3$ such that the corresponding flows $\MM^s := \MM^{g_s}$ interpolate between a round and a cylindrical singularity.
For $s \in [0,\frac12)$ the flow $\MM^s$ can be described in terms of a conventional, non-singular Ricci flow $(g^s_t)$ on $M$ and the continuity statement in Theorem~\ref{Thm_cont_dep} is equivalent to continuous dependence of this flow on $s$.
Likewise, the flows $\MM^s$ restricted to $[0,T_{\sing})$ can again be described by a continuous family of conventional Ricci flows.
The question is now what happens at the critical parameter $s = \frac12$, where the type of singularity changes.
The uniqueness property guarantees that the flows $\MM^s$ for $s \nearrow \frac12$ and $s \searrow \frac12$ must limit to the same flow $\MM^{\frac12}$.
The convergence is locally smooth, but the topology of the spacetime manifold $\MM^s$ may still change.

\subsection{Topological Applications}
Theorem~\ref{Thm_cont_dep} provides us a tool to deduce the first topological applications of Ricci flow since Perelman's work. 

The first example of such an application concerns the space of metrics of positive scalar curvature $\Met_{\PSC}(M) \subset \Met(M)$ on a manifold $M$, which is a subset of the space of all Riemannian metrics on $M$ (both spaces are equipped with the $C^\infty$-topology).
%
%
Since the positive scalar curvature condition is preserved by Ricci flow, Theorem~\ref{Thm_cont_dep} roughly implies that --- modulo singularities and the associated topological changes --- Ricci flow is a ``continuous deformation retraction'' of $\Met_{\PSC}(M)$ to the space of round metrics on $M$.
This heuristic was made rigorous by Kleiner and the author \cite{BamlerKleiner-PSC} and implied:

\begin{Theorem} \label{Thm_PSC}
For any closed 3-manifold $M$ the space $\Met_{\PSC}(M)$ is either contractible or empty.
\end{Theorem}

The study of the spaces $\Met_{\PSC}(M)$ was initiated by  Hitchin in the 70s and has led to many interesting results --- based on index theory --- which show that these spaces have non-trivial topology when $M$ is high dimensional.
Theorem~\ref{Thm_PSC} provides first examples of manifolds of dimension $\geq 3$ for which the homotopy type of $\Met_{\PSC}(M)$ is completely understood; see also prior work by Marques \cite{Marques_2012}.

A second topological application concerns the diffeomorphism group $\Diff(M)$ of a manifold $M$, i.e. the space of all diffeomorphisms $\phi : M \to M$ (again equipped with the $C^\infty$-topology).
The study of these spaces was initiated by Smale, who showed that $\Diff(S^2)$ is homotopy equivalent to the orthogonal group $O(3)$, i.e. the set of isometries of $S^2$.
More generally, we can fix an arbitrary closed manifold $M$, pick a Riemannian metric $g$ and consider the natural injection of the isometry group
\begin{equation} \label{eq_Isom_2_Diff}
 \Isom(M,g) \longrightarrow \Diff(M). 
\end{equation}
The following conjecture allows us to understand the homotopy type of $\Diff(M)$ for many important 3-manifolds.

\begin{Conjecture}[Generalized Smale Conjecture] \label{GSC}
Suppose that $(M^3,g)$ is closed and has constant curvature $K \equiv \pm 1$.
Then (\ref{eq_Isom_2_Diff}) is a homotopy equivalence.
\end{Conjecture}

This conjecture has had a long history and many interesting special cases were established using topological methods, including the case $M = S^3$ by Hatcher and the hyperbolic case by Gabai.
For more background see the first chapter of \cite{Rubinstein_book}.


An equivalent version of Conjecture~\ref{GSC} is that Theorem~\ref{Thm_PSC} remains true if we replace $\Met_{\PSC}(M)$ by the space $\Met_{K \equiv \pm 1}(M)$ of constant curvature metrics.
This was verified by Kleiner and the author \cite{BamlerKleiner2017b, BamlerKleiner-PSC}, which led to:


%
%
%

\begin{Theorem} \label{Thm_GSC}
The Generalized Smale Conjecture is true.
\end{Theorem}

The proof of Theorem~\ref{Thm_GSC} provides a unified treatment of all possible topological cases and it can also be extended to other manifolds $M$.
In addition it is independent of Hatcher's proof, so it gives an alternative proof in the $S^3$-case.

\section{Dimensions $n \geq 4$} \label{sec_dim_n_geq_4}
For a long time, most of the known results of Ricci flows in higher dimensions concerned special cases, such as K\"ahler-Ricci flows or flows that satisfy certain preserved curvature conditions.
\emph{General} flows, on the other hand, were relatively poorly understood.
Recently, however, there has been some movement on this topic --- in part, thanks to a slightly different geometric perspective on Ricci flows \cite{Bamler_HK_RF_partial_regularity,Bamler_HK_entropy_estimates,Bamler_RF_compactness}.
The goal of this section is to convey some of these new ideas and to provide an outlook on possible geometric and topological applications.

\subsection{Gradient shrinking solitons}
Let us start with the following basic question:
What would be reasonable singularity models in higher dimensions?
One important class of such models are gradient shrinking solitons (GSS).
The GSS equation concerns Riemannian manifolds $(M,g)$ equipped with a potential function $f \in C^\infty (M)$ and reads:
\[ \Ric + \nabla^2 f - \frac12 g = 0. \]
This generalization of the Einstein equation is interesting, because it gives rise to an associated selfsimilar Ricci flow:
\[ g(t) := |t| \phi_t^* g, \quad t < 0, \]
where $(\phi_t : M \to M)_{t < 0}$ is the flow of the vector field $|t| \nabla f$.

A basic example of a GSS is an Einstein metric ($\Ric = \frac12 g$), for example a round sphere.
In this case $g(t)$ just evolves by rescaling and becomes singular at time $0$.
A more interesting class of examples are round cylinders $S^{k \geq 2} \times \IR^{n-k}$, where
\[ g = 2(k-1) g_{S^k} + g_{\IR^{n-k}}, \quad 
f = \frac14 \sum_{i=k+1}^n x_i^2. \]
In this case $|t| \nabla f$ generates a family of dilations on the $\IR^{n-k}$ factor and
\[ g(t) = 2(k-1) |t| g_{S^k} + g_{\IR^{n-k}}, \]
which is isometric to $|t| g$.
In dimensions $n \leq 3$, all (non-trivial\footnote{Euclidean space $\IR^n$ equipped with $f = \frac14 r^2$ is called a trivial GSS.}) GSS are quotients of round spheres or cylinders.
However, more complicated GSS exist in dimensions $n \geq 4$.

By construction, GSS (or their associated flows, to be precise) are invariant under parabolic rescaling.
So the blow-up singularity model of the singularity at time $0$ (taken along an appropriately chosen sequence of basepoints) is equal to the flow itself.
Therefore every GSS does indeed occur as a singularity model, at least of its own flow.

Vice versa, the following conjecture, which I will keep vague for now, predicts that the converse should also be true in a certain sense.

\begin{Conjecture} \label{Conj_Folklore}
For any Ricci flow ``most'' singularity models are gradient shrinking solitons.
\end{Conjecture}

This conjecture has been implicit in Hamilton's work from the 90s and a similar result is known to be true for mean curvature flow.
In the remainder of this section I will present a resolution of a version of this conjecture.

\subsection{Examples of singularity formation} \label{subsec_Appleton}
Let us first discuss an example in order to adjust our expectations in regards to Conjecture~\ref{Conj_Folklore}.
In \cite{Appleton:2019ws}, Appleton constructs a class of 4-dimensional Ricci flows\footnote{The flows are defined on non-compact manifolds, but the geometry at infinity is well controlled.}
that develop a singularity in finite time, which can be studied via the blow-up technique from Subsection~\ref{subsec_blowup} --- this time we even allow the rescaling factors to be \emph{any} sequence of numbers $\lambda_i \to \infty$, not just $\lambda_i = |{\Rm}|^{1/2}(x_i,t_i)$.
Appleton obtains the following classification of all non-trivial blow-up singularity models:
%
%
%
\begin{enumerate}
\item The Eguchi-Hanson metric, which is Ricci flat and asymptotic to the flat cone $\IR^4/\IZ_2$.
\item The flat cone $\IR^4 / \IZ_2$, which has an isolated orbifold singularity at the origin.
\item The quotient $M_{\Bry} / \IZ_2$ of the Bryant soliton, which also has an isolated orbifold singularity at its tip.
\item The cylinder $\IR P^3 \times \IR$.
\end{enumerate}
Here the models 1., 2. \emph{have} to occur as singularity models, and it is unknown whether the models 3., 4. actually do show up.
The only gradient shrinking solitons in this list are 2., 4.
Note that the flow on $\IR^4 / \IZ_2$ is constant, but each time-slice is a  metric cone, and therefore invariant under rescaling.
So we may also view this model as a (degenerate) gradient shrinking soliton (in this case $f = \frac14 r^2$).

It is conceivable that there are Ricci flow singularities whose only blow-up models are of type 1., 2.
In addition, there are further examples in higher dimensions \cite{Stolarski:2019uk} whose only blow-up models that are gradient shrinking solitons must be singular and possibly degenerate.
This motivates the following revision of Conjecture~\ref{Conj_Folklore}.

\begin{Conjecture} \label{Conj_Folklore_revised}
For any Ricci flow ``most'' singularity models are gradient shrinking solitons that may be degenerate and may have a singular set of codimension $\geq 4$.
\end{Conjecture}

\subsection{A compactness and partial regularity theory}
The previous example taught us that in higher dimensions it becomes necessary to consider \emph{non-smooth} blow-up limits.
The usual convergence and compactness theory of Ricci flows due to Hamilton from Subsection~\ref{subsec_blowup}, which relies on curvature bounds and only produces smooth limits, becomes too restrictive for such purposes.
Instead, we need a fundamentally new \emph{compactness and partial regularity theory} for Ricci flows, which will enable us to take limits of arbitrary Ricci flows and study their structural properties.
This theory was recently found by the author and will lie at the heart of a resolution of Conjecture~\ref{Conj_Folklore_revised}.

Compactness and partial regularity theories are an important feature in many subfields of geometric analysis.
To gain a better sense for such theories, let us first review the compactness and partial regularity theory for Einstein metrics\footnote{This theory is due to Cheeger, Colding, Gromov, Naber and Tian.} --- an important special case, since every Einstein metric corresponds to a Ricci flow.
Consider a sequence $(M_i, g_i, x_i)$ of pointed, complete $n$-dimensional Einstein manifolds, $\Ric_{g_i} = \lambda_i g_i$, $|\lambda_i| \leq 1$.
Then, after passing to a subsequence, these manifolds (or their metric length spaces, to be precise) converge in the Gromov-Hausdorff sense to a pointed metric space
\[ (M_i, d_{g_i}, x_i) \xrightarrow[i \to \infty]{GH} (X, d, x_\infty). \]
If we now impose the following \emph{non-collapsing condition:}
\begin{equation} \label{eq_noncollapsing_Einstein}
 \vol B(x,r) \geq v > 0, 
\end{equation}
for some uniform $r, v > 0$, then the limit admits is a regular-singular decomposition
\[ X = \RR \,\, \dotcup \,\, \SS \]
such that the following holds:
\begin{enumerate}
\item $\RR$ can be equipped with the structure of a Riemannian Einstein manifold $(\RR, g_\infty)$ in such a way that the restriction $d|_{\RR}$ equals the length metric $d_{g_\infty}$.
In other words, $(X,d)$ is isometric to the metric completion of $(\RR, d_{g_\infty})$.
\item We have the following estimate on the Minkowski dimension of the singular set: 
\[ \dim_{\MM} \SS \leq n-4 \]
\item Every tangent cone (i.e. blow-up of $(X,d)$ pointed at the same point) is in fact a metric cone.
\item We have a filtration of the singular set
\[ \SS^0 \subset \SS^1 \subset \ldots \subset \SS^{n-4} = \SS \]
such that $\dim_{\mathcal{H}} \SS^k \leq k$ and such that every $x \in \SS^k \setminus \SS^{k-1}$ has a tangent cone that splits off an $\IR^k$-factor.
\end{enumerate}

Interestingly, compactness and partial regularity theories for other geometric equations (e.g. minimal surfaces, harmonic maps, mean curvature flow, \ldots) take a similar form.
The reason for this is that these theories rely on only  a few basic ingredients (e.g. a monotonicity formula, an almost cone rigidity theorem and an $\eps$-regularity theorem), which can be verified in each setting.
A theory for Ricci flows, however, did not exist for a long time, because these basic ingredients are --- at least \emph{a priori} --- not available for Ricci flows.
So this setting required a different approach.

Let us now discuss the new compactness and partial regularity theory for Ricci flows.
Before we begin, note that there is an additional complication:
Parabolic versions of notions like ``metric space'', ``Gromov-Hausdorff convergence'', etc. didn't exist until recently, so they --- and a theory surrounding them --- first had to be developed.
I will discuss these new notions in more detail in  Subsection~\ref{subsec_met_flow}.
For now, let us try to get by with some more vague explanations.

The first result is a compactness theorem.
Consider a sequence of pointed, $n$-dimensional Ricci flows 
\[ (M_i, (g_i(t))_{t \in (-T_i,0]},(x_i,0)), \]
where we imagine the basepoints $(x_i,0)$ to live in the final time-slices, and suppose that $T_\infty := \lim_{i \to \infty} T_i > 0$.
Then we have:

\begin{Theorem} \label{Thm_compactness_RF}
After passing to a subsequence, these flows $\IF$-converge to a pointed metric flow:
\[ (M_i, (g_i(t)), (x_i,0)) \xrightarrow[i \to \infty]{\mathbb{F}} (\XX, (\nu_{x_\infty;t})). \]
\end{Theorem}

Here a ``metric flow'' can be thought of as a parabolic version of a ``metric space''.
It is some sort of Ricci flow spacetime (as in Definition~\ref{Def_RF_spacetime}) that is allowed to have singular points; think, for example, of isolated orbifold singularities in every time-slice as in Appleton's example (see Subsection~\ref{subsec_Appleton}), or a singular point where a round shrinking sphere goes extinct.
The term ``$\mathbb{F}$-convergence'' can be thought of as a parabolic version of ``Gromov-Hausdorff convergence''.

The next theorem concerns the partial regularity of the limit 
in the non-collapsed case, which we define as the case in which
\begin{equation} \label{eq_NN_lower_bound}
 \NN_{x_i,0}(r_*^2) \geq -Y_*, 
\end{equation}
for some uniform constants $r_* > 0$ and $Y_* < \infty$.
This condition is the parabolic analogue to (\ref{eq_noncollapsing_Einstein}).
The quantity $\NN_{x,t}(r^2)$ is the pointed Nash-entropy, which is related to Perelman's $\mathcal{W}$-functional; it was rediscovered by work of Hein and Naber.
%
%
Assuming (\ref{eq_NN_lower_bound}), we have:

\begin{Theorem} \label{Thm_part_reg_RF}
There is a regular-singular decomposition
\[ \XX = \RR \, \dotcup \, \SS \]
such that:
\begin{enumerate}
\item The flow on $\RR$ can be described by a smooth Ricci flow spacetime structure.
Moreover, the entire flow $\XX$ is uniquely determined by this structure.
\item We have the following dimensional estimate
\[ \dim_{\MM^*} \SS \leq (n+2)-4. \]
\item Tangent flows (i.e. blow-ups based at a fixed point of $\XX$) are (possibly singular) gradient shrinking solitons.
\item There is a filtration $\SS^0 \subset \ldots \subset \SS^{n-2}= \SS$ with similar properties as in the Einstein case.
\end{enumerate}
\end{Theorem}

A few comments are in order here.
First, note that the fact that $\XX$ is uniquely determined by the smooth Ricci flow spacetime structure on $\RR$ is comparable to what we have observed in dimension 3 (see Subsection~\ref{subsec_RF_through_sing}), where we didn't even \emph{consider} the entire flow $\XX$.
Second, Property 2. involves a parabolic version of the Minkowski dimension that is natural for Ricci flows; a precise definition would be beyond the scope of this article.
Note that the time direction accounts for 2 dimensions, which is natural.
An interesting case is dimension $n=3$, in which we obtain that the set of singular times has dimension $\leq \frac12$; this is in line with what is known in this dimension.
Lastly, in Property~3. the role of metric cones is now taken by gradient shrinking solitons; these are analogues of metric cones, as both are invariant under rescaling.

The dimensional bounds in Theorem~\ref{Thm_part_reg_RF} are optimal.
In Appleton's example, the singular set $\SS$ may consist of an isolated orbifold point in every time-slice; so its parabolic dimension is $2 = (4+2)-4$.
On the other hand, a flow on $S^2 \times T^2$ develops a singularity at a single time and collapses to the 2-torus $T^2$, which again has parabolic dimension $2$.

\subsection{Applications}
Theorems~\ref{Thm_compactness_RF} and \ref{Thm_part_reg_RF} finally enable us to study the finite-time singularity formation and long-time behavior of Ricci flows in higher dimensions.

Regarding Conjecture~\ref{Conj_Folklore_revised}, we roughly obtain:

\begin{Theorem} \label{Thm_finite_time}
Suppose that $(M, (g(t))_{t \in [0,T)})$ develops a singularity at time $T < \infty$.
Then we can extend this flow by a ``singular time-$T$-slice'' $(M_T, d_T)$ such that the tangent flows at any $(x,T) \in M_T$ are (possibly singular) gradient shrinking solitons. 
\end{Theorem}

Regarding the long-time asymptotics, we obtain the following picture, which closely resembles that in dimension~3 compare with (\ref{eq_thick_thin_3d}) in Subsection~\ref{subsec_RFsurg}:

\begin{Theorem} \label{Thm_longtime}
Suppose that $(M, (g(t))_{t \geq 0})$ is immortal.
Then for $t \gg 1$ we have a thick-thin decomposition
\[ M = M_{\thick} (t) \,\, \dotcup  \,\, M_{\thin}(t) \]
such that the flow on $M_{\thick}(t)$ converges, after rescaling, to a singular Einstein metric $(\Ric_{g_\infty} = -g_{\infty}$) and the flow on $M_{\thin}(t)$ is collapsed in the opposite sense of (\ref{eq_NN_lower_bound}).
\end{Theorem}

Theorems~\ref{Thm_finite_time} and \ref{Thm_longtime} essentially generalize Perelman's results to higher dimensions.


\subsection{Metric flows} \label{subsec_met_flow}
So what precisely is a metric flow?
To answer this question, we will imitate the process of passing from a (smooth) Riemannian manifold $(M,g)$ to its metric length space $(M, d_g)$.
Here a new perspective on the geometry of Ricci flows will be key.

So our goal will be to turn a Ricci flow $(M, (g(t))_{t \in I})$ into a synthetic object, which we call ``metric flow''.
To do this, let us first consider the spacetime $\XX := X \times I$ and the time-slices $\XX_t := X \times \{ t \}$ equipped with the length metrics $d_t := d_{g(t)}$.
It may be tempting to retain the product structure $X \times I$ on $\XX$, i.e. to record the set of worldlines $t \mapsto (x,t)$.
However, this turns out to be unnatural.
Instead, we will view the time-slices $(\XX_t, d_t)$ as separate, metric spaces, whose points may not even be in 1-1 correspondence to some given space $X$.

It remains to record some other relation between these metric spaces $(\XX_t, d_t)$.
This will be done via the conjugate heat kernel $K(x,t;y,s)$ --- an important object in the study of Ricci flows.
For fixed $(x,t) \in M \times I$ and $s < t$ this kernel satisfies the backwards conjugate\footnote{(\ref{eq_conj_HK}) is the $L^2$-conjugate of the standard (forward) heat equation and $K(\cdot,\cdot;y,t)$ is a heat kernel placed at $(y,t)$.} heat equation on a Ricci flow background:
\begin{equation} \label{eq_conj_HK}
 (- \partial_s - \triangle_{g(s)} + R_{g(s)}) K(x,t;\cdot,s) = 0,
\end{equation}
centered at $(x,t)$.
%
%
%
This  kernel has the property that for any $(x,t)$ and $s < t$
\[ \int_M K(x,t;\cdot, s) dg(s) = 1,  \]
which motivates the definition of the following probability measures:
\[ d\nu_{(x,t);s} := K(x,t;\cdot,s) dg(s), \quad \nu_{(x,t);t} = \delta_x. \]
This is the additional information that we will record.
So we define:

\begin{Definition}
A \emph{metric flow} is (essentially\footnote{This is a simplified definition.}) given by a pair
\[ \big( (\XX_t, d_t)_{t \in I}, (\nu_{x;s})_{x \in \XX_t, s < t, s \in I} \big) \]
consisting of a family of metric spaces $(\XX_t, d_t)$ and probability measures $\nu_{x;s}$ on $\XX_s$, which satisfy certain basic compatibility relations.
\end{Definition}

So given points $x \in \XX_t$, $y \in \XX_s$ at two times $s < t$, it is not possible to say whether ``$y$ corresponds to $x$''.
Instead, we only know that ``$y$ belongs to the past of $x$ with a probability density of $d\nu_{x;s}(y)$''.
This definition is surprisingly fruitful.
For example, it is possible to use the measures $\nu_{x;s}$ to define a natural topology on $\XX$ and to understand when and in what sense the geometry of time-slices $\XX_t$ depends continuously on $t$.

The concept of metric flows also allows the definition of a natural notion of geometric convergence --- $\mathbb{F}$-convergence --- which is similar to Gromov-Hausdorff convergence.
Even better, this notion can be phrased on terms of a certain $d_{\mathbb{F}}$-distance, which is similar to the Gromov-Hausdorff distance, and the Compactness Theorem~\ref{Thm_compactness_RF} can be expressed as a statement on the compactness of a certain subset of metric flow (pairs),\footnote{Strictly speaking, $\mathbb{F}$-convergence and $d_{\mathbb{F}}$-distance concern metric flow \emph{pairs}, $(\XX, (\nu_{x;t}))$, where the second entry serves as some kind of substitute of a basedpoint.} just as in the case of Gromov-Hausdorff compactness.

\begin{figure}
\centering
\includegraphics[width=65mm]{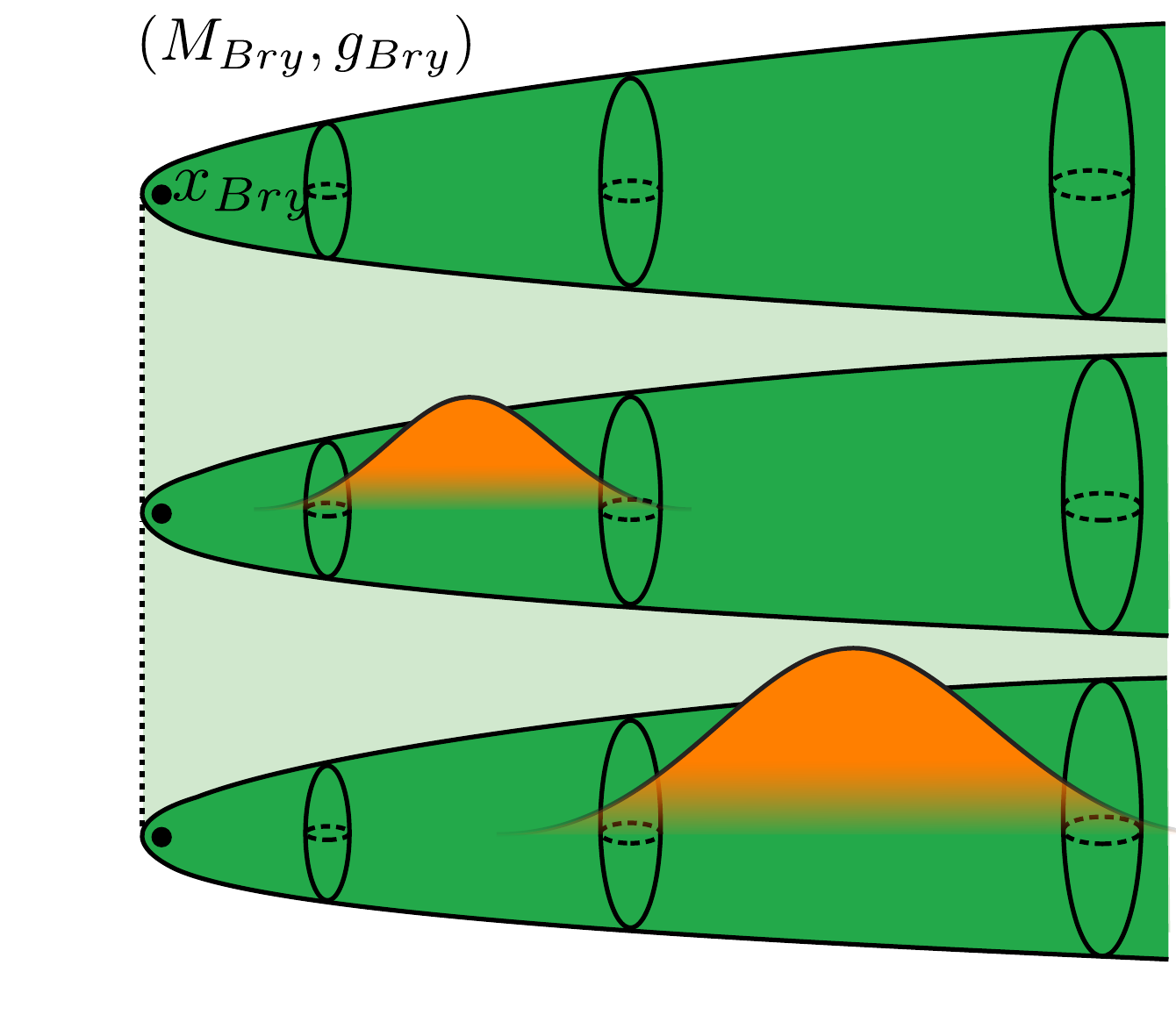}
\caption{The Bryant soliton and a conjugate heat kernel starting at $(x_{\Bry},0)$.\label{fig_Bryant_evolution}}
\end{figure}
Lastly, I will sketch an example that illustrates why it was so important that we have divorced ourselves from the concept of worldlines.
Consider the Bryant soliton $(M_{\Bry}, (g_{\Bry}(t))_{t \leq 0})$  (see Figure~\ref{fig_Bryant_evolution}) from Subsection~\ref{subsec_dumbell}.
Recall that every time-slice $(M_{\Bry}, g_{\Bry}(t))$ is isometric to the same rotationally symmetric model with center $x_{\Bry}$.
By Theorems~\ref{Thm_compactness_RF} and \ref{Thm_part_reg_RF} any pointed sequence of blow-\emph{downs} ($\lambda_i \to 0$),
\[ (M_{\Bry}, (\lambda_i^2 g_{\Bry}( \lambda_i^{-2} t))_{t \leq 0}, (x_{\Bry},0)), \]
$\IF$-converges to a pointed metric flow $\XX$ that is regular on a large set.
What is this $\IF$-limit $\XX$?
For any fixed time $t < 0$, the sequence of pointed Riemannian manifolds $(M_{\Bry}, \lambda_i^2 g_{\Bry}( \lambda_i^{-2} t), x_{\Bry})$ converges to a pointed ray of the form $([0, \infty), 0)$.
This seems to contradict Theorem~\ref{Thm_part_reg_RF}.
However, here we have implicitly used the concept of worldlines, because we have used the point $(x_{\Bry}, t)$ corresponding to the ``official'' basepoint $(x_{\Bry}, 0)$ at time $t$.
Instead, we have to focus on the ``past'' of $(x_{\Bry}, 0)$, i.e. the region of $(M_{\Bry}, \lambda_i^2 g_{\Bry}( \lambda_i^{-2} t))$ where the conjugate heat kernel $\nu_{(x_{\Bry},0);\lambda_i^{-2}t}$ is concentrated.
This region is cylindrical of scale $\sim \sqrt{|t|}$, because the conjugate heat kernel ``drifts away from the tip'' at an approximate linear rate.
In fact, one can show that the blow-down limit $\XX$ is isometric to a round shrinking cylinder that develops a singularity at time $0$.
While this may seem slightly less intuitive at first, it turns out to be a much more natural way of looking at it.

\subsection{What's next?}
This new theory demonstrates that, at least on an analytical level, Ricci flows behave similarly in higher dimension as they do in dimension 3.
However, while there are only a handful of possible singularity models in dimension 3, gaining a full understanding of all such models in higher dimensions (e.g. classifying gradient shrinking solitons) seems like an intimidating task.
Some past work in dimension 4 (e.g. by Munteanu and Wang) has demonstrated that most \emph{non-compact} gradient shrinking solitons have ends that are either cylindrical or conical.
This motivates the following conjecture:

\begin{Conjecture}
Given a closed Riemannian 4-manifold $(M,g)$ there is a certain kind of ``Ricci flow through singularities'' in which topological change occurs along cylinders or cones and in which time-slices are allowed to have isolated orbifold singularities.
\end{Conjecture}

Showing the existence of such a flow would be an analytical challenge, given that the construction in dimension 3 already filled several hundred pages.
However, I currently don't see a reason why such a flow should not exist.

Even more exciting would be the question what such a flow (or an analogue in higher dimensions) would accomplish on a topological level.
It is unlikely that it would allow us to prove the smooth Poincar\'e Conjecture in dimension 4, because even in the best possible case, such a flow seems to provide insufficient topological information.
%
A more feasible application would the $\frac{11}{8}$-Conjecture, which states that for every closed spin 4-manifold we have
\begin{equation} \label{eq_11_8}
 b_2(M) \geq \tfrac{11}8 |\sigma(M)|. 
\end{equation}
This conjecture is the missing piece in the classification of closed, simply connected, smooth 4-manifolds up to homeomorphy (due to Donaldson, Freedman and Kirby).
At least on a heuristic level, it would be suited for a Ricci flow approach since (closed) Einstein manifolds and spin gradient shrinking solitons automatically satisfy (\ref{eq_11_8}) --- due to the Hitchin-Thorpe inequality in the Einstein case or the fact that gradient shrinking solitons have positive scalar curvature.
Other potential applications would be questions concerning the topology of 4-manifolds that admit metrics of positive scalar curvature.
Lastly, there also seems to be potential in K\"ahler geometry, for example towards the Minimal Model Program and the Abundance Conjecture.

Time will tell how far Ricci flow methods will take us precisely.
Almost 20 years ago, geometric analysts and topologists were busy  digesting Perelman's work on the Poincar\'e and Geometrization Conjectures, thus closing an important chapter in the field.
Today, we have good reasons to be optimistic that further topological applications are on the horizon.
%
%
%
The field certainly still has an exciting future ahead and hopefully encourages further research and collaboration between geometers, analysts, topologists and complex geometers.

\bibliography{bibliography}{}
\bibliographystyle{amsalpha}

\end{document}